\documentclass{article}

\title 
{Axially Symmetric Generalization \\
of the Cauchy-Riemann System \\ 
and Modified Clifford Analysis }
\author{ Dmitri Bryukhov }  

\usepackage{amsthm}

\newtheorem{thm}{Theorem}[section]
\newtheorem{prop}[thm]{Proposition}

\newtheorem{cor}[thm]{Corollary}

\theoremstyle{definition}
\newtheorem{defn}[thm]{Definition}

\newtheorem{example}[thm]{Example}

\theoremstyle{remark}
\newtheorem{rem}[thm]{Remark}

\begin{document}
\maketitle

\begin{abstract}
The main aim of this paper is to describe the most adequate generalization
of the  Cauchy-Riemann system fixing properties of classical functions 
in octonionic case. An octonionic generalization of the
Laplace transform is introduced. Octonionic generalizations 
of the inversion transformation, the gamma function and 
the Riemann zeta-function are given.
\end{abstract}

$ \bf Keywords: $
{generalizations of the Cauchy-Riemann system,
functions of the octonionic variable,
octonionic Laplace transform}

\section{Introduction}

Theory of holomorphic functions $\ f=u+iv\ $ of the complex variable 
$\ z=x+iy\ $ has been developed on the basis of investigations of 
the Laplace equation on the plane $\mathbf R^2=\{(x,y)\}\ $
\begin{displaymath}
  \Delta{h} = \mathrm{div}\ \mathrm{grad}{\ h} =
 \frac{\partial{h}^2}{\partial{x}^2} + \frac{\partial{h}^2}{\partial{y}^2} = 0
\end{displaymath}
where $h\ $ - complex potential, \ and the Cauchy-Riemann system 
\begin{displaymath}
  \left\{
    \begin{array}{l}
      \frac{\partial{u}}{\partial{x}}-\frac{\partial{v}}{\partial{y}}=0 \\
      \frac{\partial{u}}{\partial{y}}=-\frac{\partial{v}}{\partial{x}}
    \end{array}
  \right.
\end{displaymath}
where $\ u(x,y)=\frac{\partial{h}}{\partial{x}},
 \ v(x,y)=-\frac{\partial{h}}{\partial{y}}\ $ (see, e.g. \cite {LavSh}).

Leutwiler \cite{Leut:CV17} considered 
the remarkable hyperbolic version of the Laplace equation
in $\mathbf R^{n+1}=\{(x_0,x_1,...,x_n)\}$ 
\begin{equation}
  x_n\Delta{h}-(n-1)\frac{\partial{h}}{\partial{x_n}}=0 \ \ \ \ \
  (\Delta = \frac{\partial{}^2}{\partial{x_0}^2}+
  \frac{\partial{}^2}{\partial{x_1}^2}+...+
  \frac{\partial{}^2}{\partial{x_n}^2}). 
\end{equation}

\begin{rem}
  It is easily seen, if $\ \ x_n\neq0\ $ then \\
\begin{displaymath}
 x_n\Delta{h}-(n-1)\frac{\partial{h}}{\partial{x_n}}=
  x_n^n\mathrm{div}\ (x_n^{1-n}\mathrm{grad}{\ h})=0.
\end{displaymath}
\end{rem} 

 The collection of (n+1) real $C^2$-functions $\ u_0=u_0(x_0,x_1,...,x_n)$, \\ 
$u_1=u_1(x_0,x_1,...,x_n),\ ...,\ u_n=u_n(x_0,x_1,...,x_n)$, \\
where  $\ u_0=\frac{\partial{h}}{\partial{x_0}},
 \ \ u_1=-\frac{\partial{h}}{\partial{x_1}},..., 
 \ \ u_n=-\frac{\partial{h}}{\partial{x_n}},\ \ $ in this case \\ 
satisfies the asymmetric system $\ (H_n)$
\begin{equation}
  \left\{
    \begin{array}{l}
      x_n(\frac{\partial{u_0}}{\partial{x_0}}-
      \frac{\partial{u_1}}{\partial{x_1}}-...-
      \frac{\partial{u_n}}{\partial{x_n}}) + (n-1)u_n=0 \\
      \frac{\partial{u_0}}{\partial{x_m}}=-\frac{\partial{u_m}}{\partial{x_0}}
      \ \ \ \ \ \ \ (m=1,...,n)  \\
      \frac{\partial{u_l}}{\partial{x_m}}=\ \ \frac{\partial{u_m}}{\partial{x_l}}
      \ \ \ \ \  (l,m = 1,...,n) 
     \end{array}
  \right.\label{eq:H_n-system}
\end{equation}
Leutwiler investigated various classes of solutions (\ref {eq:H_n-system}) 
connected with the (universal) Clifford algebra $\mathbf C_n$, 
in particular with the quaternion algebra $\mathbf {H=C_2}$.

\begin{rem}
 The Clifford algebra $\mathbf C_3$ (without division)    
and the octonion algebra $\mathbf O$ (with division) aren't equivalent.  
\end{rem}

The Laplace-Beltrami equation in $\mathbf R^3=\{(x,y,t)\}$
\begin{displaymath}
  t\Delta{h}-\frac{\partial{h}}{\partial{t}}=0 \ \ \ \ \ \
  (\Delta = \frac{\partial{}^2}{\partial{x}^2}+
  \frac{\partial{}^2}{\partial{y}^2}+
  \frac{\partial{}^2}{\partial{t}^2}) 
\end{displaymath}
was exploited in papers \cite{Leut:CV20},\cite{Leut:FM7},\cite{Leut:EM14} 
to obtain the asymmetric system $\ (H)$  
\begin{displaymath}
  \left\{
    \begin{array}{l}
      t(\frac{\partial{u}}{\partial{x}}-
      \frac{\partial{v}}{\partial{y}}-
      \frac{\partial{w}}{\partial{t}}) + w=0 \\
      \frac{\partial{u}}{\partial{y}}=-\frac{\partial{v}}{\partial{x}}, \
      \frac{\partial{u}}{\partial{t}}=-\frac{\partial{w}}{\partial{x}}, \
      \frac{\partial{v}}{\partial{t}}=\frac{\partial{w}}{\partial{y}},
    \end{array} 
  \right.
\end{displaymath}
where  $\ u=\frac{\partial{h}}{\partial{x}},
 \ \ v=-\frac{\partial{h}}{\partial{y}}, 
 \ \ w=-\frac{\partial{h}}{\partial{t}}.$  
There were introduced solutions in the form of various elementary functions
$\ f=u+iv+jw\ $ of the reduced quaternionic variable $\ z=x+iy+jt\ $.

The Laplace-Beltrami equation in $\mathbf R^4=\{(x,y,t,s)\}$
\begin{displaymath}
 s\Delta{h}-2\frac{\partial{h}}{\partial{t}}=0 \ \ \ \ \ \
  (\Delta = \frac{\partial{}^2}{\partial{x}^2}+
  \frac{\partial{}^2}{\partial{y}^2}+
  \frac{\partial{}^2}{\partial{t}^2}+
  \frac{\partial{}^2}{\partial{s}^2})
\end{displaymath}
was applied in paper \cite {HeLe:1996} to generalize the system $\ (H)$
\begin{displaymath}
  \left\{
    \begin{array}{l}
      s(\frac{\partial{u}}{\partial{x}}-
      \frac{\partial{v}}{\partial{y}}-
      \frac{\partial{w}}{\partial{t}}-
      \frac{\partial{r}}{\partial{s}}) + 2r=0 \\
      \frac{\partial{u}}{\partial{y}}=-\frac{\partial{v}}{\partial{x}}, \
      \frac{\partial{u}}{\partial{t}}=-\frac{\partial{w}}{\partial{x}}, \
      \frac{\partial{u}}{\partial{s}}=-\frac{\partial{r}}{\partial{x}} \\
      \frac{\partial{v}}{\partial{t}}=\ \frac{\partial{w}}{\partial{y}}, \
      \frac{\partial{v}}{\partial{s}}=\ \frac{\partial{r}}{\partial{y}}, \
      \frac{\partial{w}}{\partial{s}}=\ \frac{\partial{r}}{\partial{t}}
    \end{array} 
  \right.
\end{displaymath}
and to obtain solutions in the form of elementary functions
$\ f=u+iv+jw+kr\ $ of the quaternionic variable $\ z=x+iy+jt+ks\ $.

Interesting papers 
on octonion analysis \cite{LiKaPe:2001},\cite{LiPe:2001}
and on functions of the octonionic variable \cite{ScTiTo:1997}
have appeared last years. 
However generalizations of the Cauchy-Riemann system describing
properties of solutions in the form of functions of the octonionic variable
haven't been considered there. 

\section { On Axial Symmetry and Solutions Associated 
to Holomorphic Functions in $\mathbf R^{n+1}$ }

Leutwiler \cite{Leut:CV17} defined an important class 
of solutions associated to classical holomorphic functions 
of the system (\ref{eq:H_n-system}) 
(in particular $x^k,$ where $k\in \mathbf N,\ $ $e^x,\ $ $\ln x$) 
and  gave axially symmetric conditions
\begin{equation}
{x_l}{u_m}={x_m}{u_l} \ \ \ \ \ (l,m=1,...,n), 
\label{spec.cond}
\end{equation}
characterizing this class, at least locally. 

Let us introduce the following second order elliptic equation
in $\mathbf R^{n+1}=\{(x_0,x_1,...,x_n)\}$
\begin{equation}
  (x_1^2+...+x_n^2)\Delta{h}-
  (n-1)(x_1\frac{\partial{h}}{\partial{x_1}}+...+ 
        x_n\frac{\partial{h}}{\partial{x_n}})=0 
\label{eq:sigma-n}
\end{equation}
\begin{rem}
  If $\ \ (x_1^2+...+x_n^2)\neq0,\ $ then \\

  $(x_1^2+...+x_n^2)\Delta{h}-
  (n-1)(x_1\frac{\partial{h}}{\partial{x_1}}+...+ 
        x_n\frac{\partial{h}}{\partial{x_n}}) =$
 
  $(x_1^2+...+x_n^2)^{\frac{n+1}{2}}
  \mathrm{div}[(x_1^2+...+x_n^2)^{\frac{1-n}{2}}\mathrm{grad}{\ h}] = 0. \\ $
\end{rem}  

The collection of (n+1) real $C^2$-functions  $\ u_0=u_0(x_0,x_1,...,x_n),$
\\ $\ u_1=u_1(x_0,x_1,...,x_n),\ ...,
\ u_n=u_n(x_0,x_1,...,x_n)$, \\
where  $\ u_0=\frac{\partial{h}}{\partial{x_0}},
 \ \ u_1=-\frac{\partial{h}}{\partial{x_1}},..., 
 \ \ u_n=-\frac{\partial{h}}{\partial{x_n}},\ \ $ in this case \\
satisfies an axially symmetric system $\ (A_n)$ 
\begin{equation}
  \left\{
    \begin{array}{l}
      (x_1^2+...+x_n^2)(\frac{\partial{u_0}}{\partial{x_0}}-
      \frac{\partial{u_1}}{\partial{x_1}}-...-
      \frac{\partial{u_n}}{\partial{x_n}})+(n-1)(x_1u_1+...+x_nu_n)=0 \\
      \frac{\partial{u_0}}{\partial{x_m}}=-\frac{\partial{u_m}}{\partial{x_0}}
      \ \ \ \ \ \ \ (m=1,...,n)  \\
      \frac{\partial{u_l}}{\partial{x_m}}=\ \ \frac{\partial{u_m}}{\partial{x_l}}
      \ \ \ \ \  (l,m = 1,...,n) 
     \end{array}
  \right.\label{eq:A_n-system}
\end{equation}

The singular hyperplane plays an essential role in modified Clifford analysis. 
\begin{defn} 
The subspace $\mathbf R^n=\{(x_0,x_1,...,x_{n-1})\}$ of the Euclidean space
$\mathbf R^{n+1}=\{(x_0,x_1,...,x_n)\}$
is called the singular hyperplane $\mathbf {[x_n=0]}$. 
\end{defn}

\begin{thm}
In any point in ${\mathbf R^{n+1}}\setminus{\mathbf {[x_n=0]}}\ $
a collection of (n+1) real $C^2$-functions $\ (u_0,u_1,...,u_n)\ $ 
with conditions (\ref{spec.cond}) is a solution of the system 
(\ref{eq:H_n-system}) 
if and only if the collection $\ (u_0,u_1,...,u_n)\ $ 
is a solution of the system (\ref{eq:A_n-system}). 
\end{thm}

\begin{proof}
 \ \ Let $\ x_n\neq0$. \\
 If a solution $(u_0,u_1,...,u_n)$ of the system (\ref{eq:H_n-system})   
satisfies conditions (\ref{spec.cond}) then 

 $x_n(\frac{\partial{u_0}}{\partial{x_0}}-
      \frac{\partial{u_1}}{\partial{x_1}}-...-
      \frac{\partial{u_n}}{\partial{x_n}}) + (n-1)u_n =$

 $ (\sum\limits_{m=1}^n{x_m^2})
  x_n(\frac{\partial{u_0}}{\partial{x_0}}-
      \frac{\partial{u_1}}{\partial{x_1}}-...-
      \frac{\partial{u_n}}{\partial{x_n}})
  + (n-1)u_n(\sum\limits_{m=1}^n{x_m^2}) =$

 $ (\sum\limits_{m=1}^n{x_m^2})
  x_n(\frac{\partial{u_0}}{\partial{x_0}}-
      \frac{\partial{u_1}}{\partial{x_1}}-...-
      \frac{\partial{u_n}}{\partial{x_n}})
  + (n-1)x_n(\sum\limits_{m=1}^n{x_mu_m}) = 0$ \\
and we have the first equation of the system (\ref{eq:A_n-system}).   

 If a solution $(u_0,u_1,...,u_n)$ of the system (\ref{eq:A_n-system})   
satisfies conditions (\ref{spec.cond}) then 

  $ (x_1^2+...+x_n^2)(\frac{\partial{u_0}}{\partial{x_0}}-
    \frac{\partial{u_1}}{\partial{x_1}}-...-
    \frac{\partial{u_n}}{\partial{x_n}})+(n-1)(x_1u_1+...+x_nu_n) = $

  $ x_n(x_1^2+...+x_n^2)(\frac{\partial{u_0}}{\partial{x_0}}-
    \frac{\partial{u_1}}{\partial{x_1}}-...-
    \frac{\partial{u_n}}{\partial{x_n}})+(n-1)(x_1u_1+...+x_nu_n)x_n = $

  $ x_n(x_1^2+...+x_n^2)(\frac{\partial{u_0}}{\partial{x_0}}-
    \frac{\partial{u_1}}{\partial{x_1}}-...-
    \frac{\partial{u_n}}{\partial{x_n}})+(n-1)u_n(x_1^2+...+x_n^2) = 0$ \\
and we have the first equation of the system (\ref{eq:H_n-system}).   
\end{proof}

\begin{cor}
All solutions associated to classical holomorphic functions 
on the singular hyperplane $\mathbf {[x_n=0]}\ $ 
(except lower singular hyperplane 
$\mathbf {[x_{n-1}]}\equiv{\mathbf R^{n-1}}$)
are solutions of the axially symmetric equation 
\begin{displaymath}
    (x_1^2+...+x_{n-1}^2)(\frac{\partial{u_0}}{\partial{x_0}}-
    \frac{\partial{u_1}}{\partial{x_1}}-...-
    \frac{\partial{u_{n-1}}}{\partial{x_{n-1}}})
    +(n-2)(x_1u_1+...+x_{n-1}u_{n-1})=0 
     \label{eq:A-hyperplane equation}
\end{displaymath}
\end{cor}

\begin{proof}
\ \ In according with the previous theorem solutions
associated to classical holomorphic functions on 
$\mathbf {[x_n=0]}\equiv\mathbf{ R^n}$
(except the subspace $\mathbf R^{n-1}=\{(x_0,x_1,...,x_{n-2})\}$)
satisfy the system $(A_{n-1})$. The first equation of the system $(A_{n-1})$ 
coincides with the equation (\ref{eq:A-hyperplane equation}).   
\end{proof}
Note that every collection $(u_0,u_1,...,u_n)$ with conditions (\ref{spec.cond})  
has the component $u_n=0$ on the singular hyperplane $\mathbf {[x_n=0]}$. 

Thus the system (\ref{eq:A_n-system}) can be interpreted as 
a natural axially symmetric generalization of the Cauchy-Riemann system 
having solutions associated to classical holomorphic functions 
in ${\mathbf R^{n+1}}\setminus{\mathbf R}.$ 

\begin{rem}
If a collection of (n+1) real $C^2$-functions $\ (u_0,u_1,...,u_n)\ $ 
is a solution of the system (\ref{eq:A_n-system}) in $\mathbf R^{n+1},\ $
then for every $\mathbf R^2=\{(x_0,x_m)\}\ (m=1,...,n)\ $
in $\mathbf R^2\setminus{\mathbf R}\ $ 
the solution $\ (u_0,u_1,...,u_n)\ $ satisfies the simple relation  
\begin{displaymath}
    x_m(\frac{\partial{u_0}}{\partial{x_0}}-
    \frac{\partial{u_1}}{\partial{x_1}}-...-
    \frac{\partial{u_n}}{\partial{x_n}})+(n-1)u_m = 0 
     \label{eq:R^2-equation}
\end{displaymath}
\end{rem}

\section{ On Real-Valued Originals and the Octonionic \\
Generalization of the Laplace Transform} 
Recall \cite{LiKaPe:2001} that the octonion algebra $\mathbf O$
is an alternative, non-associative division algebra over $\mathbf R$
with $e_0=1$ and the basic octonion units $e_1,e_2,e_3,e_4,e_5,e_6,e_7$,
where $e_3=e_1e_2,\ e_5=e_1e_4,\ e_6=e_2e_4,\ e_7=e_3e_4.$
Thus 

$x=x_0+\sum\limits_{m=1}^7{x_me_m} =
x_0+x_1e_1+x_2e_2+x_3e_3 + (x_4+x_5e_1+x_6e_2+x_7e_3)e_4$.

 If $x\notin\mathbf R$ we can use the polar form 
\begin{equation}
x=x_0 +\sum\limits_{m=1}^7{x_me_m} =
|x|(\cos{\varphi}+I(x)\sin{\varphi})=|x|e^{I(x)\varphi},
\end{equation}
 where $\ \ \ \ \ I(x)= 
\frac{x_1e_1+x_2e_2+x_3e_3+x_4e_4+x_5e_5+x_6e_6+x_7e_7}
{\sqrt{x_1^2+x_2^2+x_3^2+x_4^2+x_5^2+x_6^2+x_7^2}}\ \ \ \ \
(I(x)^2=-1)$,

$\ \ \ \ \ \ \ \ \ \ \varphi=\arccos\frac{x_0}
{\sqrt{x_1^2+x_2^2+x_3^2+x_4^2+x_5^2+x_6^2+x_7^2}}\ \ \ \ \ \ \ \ \ \
(0<\varphi<\pi)$.

Then for any $x\notin\mathbf R$ 
\begin{equation}
\ln{x}=\ln|x|+I(x)\varphi\ \ \ \ (principal\ value)
\end{equation}
and for any $n\in{\mathbf N}$
\begin{equation}
x^n=|x|^n(\cos{n\varphi}+I(x)\sin{n\varphi}).
\end{equation}

Analogously \cite{Leut:CV17} formula
\begin{displaymath}
e^{I(x)\lambda}=\cos\lambda+I(x)\sin\lambda,
\end{displaymath}
where $\ \ \lambda\in{\mathbf R}\ $ 
(the octonionic analog of Euler's relation), 
\\ in case of 
$\ \lambda=\rho={\sqrt{x_1^2+x_2^2+x_3^2+x_4^2+x_5^2+x_6^2+x_7^2}} $ \\ and
$\ \ \ \ \ \ \ \ I(x)\rho=x_1e_1+x_2e_2+x_3e_3+x_4e_4+x_5e_5+x_6e_6+x_7e_7,$ 
\\ has as a consequence the following formula 
\begin{equation}
e^x=e^{x_0}e^{I(x)\rho}=e^{x_0}(\cos{\rho}+I(x)\sin{\rho}). 
\end{equation}
The octonionic inversion is described by the simple relation
\begin{equation}
x^{-1}=\frac{\bar{x}}{|x|^2}=
\frac{x_0-\sum\limits_{m=1}^7{x_me_m}}{|x|^2}=
|x|(\cos{\varphi}-I(x)\sin{\varphi})=|x|e^{-I(x)\varphi}.
\end{equation}
Then
\begin{equation}
x^{-n}=|x|^{-n}(\cos{n\varphi}-I(x)\sin{n\varphi})
=|x|^{-n}e^{-I(x)n\varphi}.
\end{equation}

\begin{rem}
As is easily seen, elementary functions 
$x^n, \ln{x}, e^x, x^{-n}\ $ of the octonionic variable $x\ $
satisfy the symmetric conditions
$\ \ {u_l}{x_m}={u_m}{x_l}\ \ \ \ (l,m=1,...,7)$.
Therefore for $m=1,...,7\ $ a condition $\ x_m=0\ $ implies $\ u_m=0$.  
\end{rem}

It is directly verified that elementary functions 
$x^n, \ln{x}, e^x, x^{-n}\ $ of the octonionic variable $x\ $
generate solutions of the system $(A_7)$. 
Moreover the system $(A_7)$ is linear therefore any linear combinations 
of these elementary functions generate solutions too.

\begin{defn}
  A real-valued function $\eta(\tau)$ of a real variable $\tau$
  is called an original, if
  \begin{enumerate}
  \item  $\eta(\tau)$ complies with the H\"{o}lder's condition for every
    $\tau$ except some points $\tau =
    \tau^1_{\eta},\tau^2_{\eta},\ldots$ (there exists a finite
    quantity or zero of such points for every finite interval), where
    the function $\eta(\tau)$ has gaps of the first kind,
  \item $\eta(\tau) = 0$ for all $\tau<0$,
  \item there exist constants $B_{\eta}>0,x^0_{\eta}\geq0:\ $ for all
    $\tau$ $\ |\eta(\tau)| < B_{\eta} e^{x^0_{\eta}\tau}$.
  \end{enumerate}

 The H\"{o}lder's condition  for the  function $\eta(\tau)$
has the form:
for every $\ \tau,\ $ there exist constants $\ A_{\eta}>0,\
0<\lambda_{\eta}\leq1, \ \delta_{\eta}>0\ $ so that
$|\eta(\tau+\delta)-\eta(\tau)|
\leq A_{\eta}|\delta|^{\lambda_{\eta}}$
for every $\delta,\ |\delta|\leq\delta_{\eta}$.
\end{defn}

\begin{rem}
  It is well known that in complex case for every original
  $\eta(\tau)$ the Laplace transform exists in the area $Re\
  z=x>x^0_{\eta}$. Similar property plays an important role in 
  octonionic case too.
\end{rem}

\begin{defn}
  For every original $\eta(\tau)\ $ a function of an octonionic
  variable
\begin{equation}
 \mathcal{L}[\eta(\tau)](x) =\int_{0}^{\infty}
  e^{-x\tau}\eta(\tau) d\tau 
\label{def:LaTr}  \end{equation}
is called an octonionic generalization of the Laplace transform 
(or simply the octonionic Laplace transform).
\end{defn}
\begin{rem}
  It is clear that $\ \mathcal{L}[\eta(\tau)](x)$ 
  $= \int\limits_{0}^{\infty}
  e^{-x\tau}\eta(\tau) d\tau$ $= \int\limits_{-\infty}^{+\infty}
  e^{-x\tau}\eta(\tau) d\tau$.
\end{rem}

\begin{prop}
  The octonionic Laplace transform $\ \mathcal{L}[\eta(\tau)](z)$ for
  every real original $\ \eta(\tau) $ defines a solution associated to
a classical holomorphic function.
\end{prop}

\begin{proof}
  \ \ Let $\ \mathcal{L}[\eta(\tau)](z)$
  $\ =u_0+\sum\limits_{m=1}^7{u_me_m}$. \\
  The octonionic exponential function defines a solution $(u_0,u_1,...,u_7)\ $
associated to the classical exponential function of the system $(A_7)$.

 Besides, \\

  $\frac{\partial{}}{\partial{x_m}}
  \int\limits_{0}^{\infty}e^{-x\tau}\eta(\tau) d\tau$
  $=\int\limits_{0}^{\infty}\frac{\partial{}}{\partial{x_m}}
  e^{-x\tau}\eta(\tau) d\tau \ \ \ \ \ (m=0,1,...,7).\\ $ 

  Then we can directly calculate that \\

   $(x_1^2+...+x_7^2)(\frac{\partial{u_0}}{\partial{x_0}}-
    \frac{\partial{u_1}}{\partial{x_1}}-...-
    \frac{\partial{u_7}}{\partial{x_7}})+6(x_1u_1+...+x_7u_7)=0$, 

   $\frac{\partial{u_0}}{\partial{x_m}}=-\frac{\partial{u_m}}{\partial{x_0}}
    \ \ \ \ \ (m=1,...,7)$,  

   $\frac{\partial{u_l}}{\partial{x_m}}=\ \ \frac{\partial{u_m}}{\partial{x_l}}
    \ \ \  (l,m = 1,...,7). $
\end{proof}

\begin{example}
  The original $\ \eta(\tau)= \left\{\begin{array}{l}$ $ 1, \
      \tau\geq0 \\ $ $ 0, \ \tau<0 $ $\end{array} \right. $
  \ implies  \\
  $\mathcal{L}[\eta(\tau)](x)=x^{-1}.$
\end{example}
\begin{example}
  The
  original $\ \eta(\tau)= \left\{\begin{array}{l}$ $ cos\ \omega\tau,
      \tau\geq0 \\ $ $ 0, \ \tau<0 $ $\end{array} \right.$
  \ implies  \\
  $\mathcal{L}[\eta(\tau)](x)=x(x^2+\omega^2)^{-1}.$
\end{example}
\begin{example}
  The original $\ \eta(\tau)= \left\{\begin{array}{l}$ $
      \tau^a, \tau\geq0 \\ $ $ 0, \ \tau<0 $ $\end{array} \right. $ \
  for every $\ a>0$ \\ implies  $\
  \mathcal{L}[\eta(\tau)](x)=\Gamma(a+1)x^{-a-1},$ \\
where $\Gamma(a+1)$ denotes the classical gamma function of a real argument.
\end{example}

\begin{rem}
Examples aren't correct for Clifford algebras $\mathbf C_n\ (n \geq 3)$.
\end{rem}

 It isn't difficult to introduce an octonionic generalization of 
the two-sided (or bilateral) Laplace transform (see, e.g. \cite{PolBr}) 
for real originals, if $\ \ \eta(\tau)\neq 0\ $ $(\tau < 0).$ 
\begin{defn}
  For every original $\eta(\tau)\ $ a function of an octonionic
  variable
\begin{equation}
 \mathcal{L}^{(2)}[\eta(\tau)](x) =\int_{-\infty}^{\infty}
  e^{-x\tau}\eta(\tau) d\tau 
\label{def:LaTr2}  \end{equation}
is called an octonionic generalization of the two-sided Laplace transform 
(or simply the octonionic two-sided Laplace transform).
\end{defn}

Thus, natural generalizations of many classical functions
(see, e.g. \cite{Titch:R},\cite{WW}) can be obtained.

\begin{example}
An octonionic generalization of the gamma function. 

Let $\ \ x=x_0+\sum\limits_{m=1}^7{x_me_m}, \ \ x_0>0.$
 \begin{eqnarray*}
  \int\limits_{-\infty}^{+\infty} e^{-x\tau}e^{-e^{\tau}}d\tau
  = \int\limits_{0}^{\infty} {\tau_1}^{-x-1}e^{-\tau_1}d\tau_1, 
  \ \ \ where \ \tau_1=e^{\tau}, \ d\tau_1=e^{\tau}d\tau. 
 \end{eqnarray*}

 Then we can correctly define  
  \begin{equation}
    \Gamma(x)= \int\limits_{0}^{\infty}
    {\tau_1}^{x-1}e^{-\tau_1}d\tau_1=
    \int\limits_{-\infty}^{+\infty} e^{x\tau}e^{-e^{\tau}}d\tau.
  \end{equation}
\end{example}

\begin{example}
An octonionic generalization of the Riemann zeta-function.

Let $\ \ x=x_0+\sum\limits_{m=1}^7{x_me_m}, \ \ x_0>1.$
  \begin{eqnarray*}
    \lefteqn{\ \ \ \ \ \ \ \ \ \ \int\limits_{-\infty}^{+\infty}
      {\frac{e^{-x\tau}d\tau}{e^{e^{\tau}}-1}}=
    \int\limits_{0}^{\infty}
      {\frac{{\tau_1}^{-x-1}d\tau_1}{e^{\tau_1}-1}} }\\
    & = &\int\limits_{0}^{\infty}
    {\tau_1}^{-x-1}(\sum_{n=1}^{\infty}e^{-n\tau_1})d\tau_1
     = \sum\limits_{n=1}^{\infty}
    \int\limits_0^{\infty}{\tau_1}^{-x-1}e^{-n\tau_1}d\tau_1 \\
    & = &\sum\limits_{n=1}^{\infty}(n^{x}
    \int\limits_0^{\infty}{\tau_2}^{-x-1}e^{-\tau_2})d\tau_2
     = (\sum\limits_{n=1}^{\infty}n^{x})
    \int\limits_0^{\infty}{\tau_2}^{-x-1}e^{-\tau_2}d\tau_2 \\
    & = & \int\limits_0^{\infty}{\tau_2}^{-x-1}e^{-\tau_2}d\tau_2 \
    (\sum\limits_{n=1}^{\infty}n^{x}),
  \ \ \ where \ \tau_1=e^{\tau},\ d\tau_1=e^{\tau}d\tau,\ \tau_2=n\tau_1.
  \end{eqnarray*}

  Then we can correctly define
  \begin{equation}
    \zeta(x)= \sum\limits_{n=1}^{\infty}n^{-x}=\Gamma^{-1}(x)
    \int\limits_0^{\infty}
    {\frac{{\tau_1}^{x-1}d\tau_1}{e^{\tau_1}-1}}.
  \end{equation}
\end{example}

\section {Boundary Value Problems and \\ 
Functions of the Octonionic Variable}

  Second order elliptic equations in divergence form 
 have various interesting applications in mathematical physics
(see, e.g. \cite{Landis}).
For a stationary temperature field $h\ $ the function 
$\overline{f}=\mathrm{grad}\ h\ $ could be interpreted as
the temperature gradient in $\mathbf R^{n+1}$. 
 If $\chi$ is the coefficient of heat conductivity, then the
 equation of heat conduction has the form 
$\ \ \mathrm{div}({\chi}\ \mathrm{grad}{\ h})=0.$

 The equation with an asymmetric $\ \chi$-distribution 
\begin{equation}
   \mathrm{div}(x_n^{1-n}\mathrm{grad}{\ h})=0  
\end{equation}
is equivalent to the system (\ref{eq:H_n-system}), at least  
in simply connected domains $\ \Lambda\subset\Omega$
$(\Lambda\subset{\mathbf R^{n+1}},\ x_n\neq0).$

 The equation with an axially symmetric $\ \chi$-distribution 
\begin{equation}
 \mathrm{div}[(x_1^2+...+x_n^2)^{\frac{1-n}{2}}\mathrm{grad}{\ h}] = 0
\label {eq:symm}
\end{equation}
is equivalent to the system (\ref{eq:A_n-system}), at least  
in simply connected domains $\ \Lambda\subset\Omega$
$(\Lambda\subset{\mathbf R^{n+1}},\ x_1^2+...+x_n^2\neq0).$

For example, the classic function of the octonionic variable
$\overline{f(x)}=\overline{x^{-1}}=\mathrm{grad}{\ h},\ x\neq0,$ 
is the inversion transformation in $\mathbf R^8$ (see, e.g. \cite {GilTru}). 
It could be interpreted as a nonstandard generalization of the plane field 
of the unit source in case of varying coefficient of heat conductivity
$\ \chi$.

\begin{thm}[Uniqueness]
  Assume that a simply connected domain $\ \Lambda\subset{\mathbf R^{n+1}}$
 $(\Lambda\cap\mathbf R=\emptyset)\ $ 
 has the  $C^2$-boundary $\partial{\Lambda}$. 
  Let $\ P=(P_0,P_1,...,P_n),\ |P|=1,\  $ 
 is outer unit normal to $\partial{\Lambda}$.
 Assume that there exist two functions
  $\ \hat f=\hat f(x)=\hat u_0+e_1\hat u_1+...+e_n\hat u_n\ $ and
  $\ \check f=\check f(x)=\check u_0+e_1\check u_1+...+e_n\check u_n\ $ 
determining regular in $\Lambda$ solutions of 
the first boundary value problem for the system (\ref{eq:A_n-system})
  \begin{displaymath}
    u_0|_{\partial{\Lambda}}=\psi_0,\
    u_1|_{\partial{\Lambda}}=-\psi_1,\ ...,\ 
    u_n|_{\partial{\Lambda}}=-\psi_n,\ \
    \psi=(\psi_0,\psi_1,...,\psi_n)
\in{C^0}(\partial{\Lambda}). \\
  \end{displaymath}

  If there doesn't exist a point $\ x^0\in{\partial{\Lambda}},\ $ 
  where $(P,\ \psi) = \sum\limits_{m=0}^n{P_m\psi_m}=0,\ $
  then $\ \hat f=\check f$.
\end{thm}

\begin{proof}
  The first boundary value problem
  \begin{displaymath}
    u_0|_{\partial{\Lambda}}=\psi_0,\
    u_1|_{\partial{\Lambda}}=-\psi_1,\ ...,\ 
    u_n|_{\partial{\Lambda}}=-\psi_n
  \end{displaymath}
  for the system (\ref{eq:A_n-system})
is equivalent to the third boundary value problem
  \begin{displaymath}
    \frac{\partial{h}}{\partial{x_0}}|_{\partial{\Lambda}}=\psi_0,\
    \frac{\partial{h}}{\partial{x_1}}|_{\partial{\Lambda}}=\psi_1,\ ...,\
    \frac{\partial{h}}{\partial{x_n}}|_{\partial{\Lambda}}=\psi_n
  \end{displaymath}
  for the equation (\ref{eq:symm}).
  Let us have
  \begin{eqnarray*}
    \hat u_0=\frac{\partial{\hat h}}{\partial{x_0}}, \quad
    \hat u_1=-\frac{\partial{\hat h}}{\partial{x_1}},\ ..., \quad
     \hat u_n=-\frac{\partial{\hat h}}{\partial{x_n}}, \\
    \check u_0=\frac{\partial{\check h}}{\partial{x_0}}, \quad
    \check u_1=-\frac{\partial{\check h}}{\partial{x_1}},\ ..., \quad
     \check u_n=-\frac{\partial{\check h}}{\partial{x_n}}. 
  \end{eqnarray*}
  Then for the function $\ h = \hat h - \check h\ $ we will obtain
  \begin{eqnarray*}
   \hat u_0 - \check u_0=\frac{\partial{h}}{\partial{x_0}},&&
   \hat u_1 - \check u_1=-\frac{\partial{h}}{\partial{x_1}},\ ..., \quad
   \hat u_n - \check u_n=-\frac{\partial{h}}{\partial{x_n}}\ \qquad \textrm{and}
   \\
   \frac{\partial{h}}{\partial{x_0}}|_{\partial{\Lambda}}=0,&&
    \frac{\partial{h}}{\partial{x_1}}|_{\partial{\Lambda}}=0,\ ...,\quad
    \frac{\partial{h}}{\partial{x_n}}|_{\partial{\Lambda}}=0.
  \end{eqnarray*}
  If there doesn't exist a point $\ x^0\in{\partial{\Lambda}},\ $ 
  where $(P,\ \psi) = 0,\ $
  then the homogeneous boundary value problem
  \begin{displaymath}
    \frac{\partial{h}}{\partial{x_0}}|_{\partial{\Lambda}}=0,\quad
    \frac{\partial{h}}{\partial{x_1}}|_{\partial{\Lambda}}=0,\ ..., \quad
    \frac{\partial{h}}{\partial{x_n}}|_{\partial{\Lambda}}=0
  \end{displaymath}
  for the equation (\ref{eq:symm}) can only have a constant solution 
(see, e.g. \cite{Bit:BoEl}).
  Hence $\ h=const\ $ and 
  $\ \hat u_0 - \check u_0=0,
   \ \hat u_1 - \check u_1=0,\ ...,\
   \ \hat u_n - \check u_n=0.$
\end{proof}

Thus under suitable conditions classic functions of the 
octonionic variable can define regular solutions $\ h=h(x)$ 
of the third boundary value problem for the equation 
\begin{equation}
 \mathrm{div}[(x_1^2+x_2^2+x_3^2+x_4^2+x_5^2+x_6^2+x_7^2)
^{-3}\mathrm{grad}{\ h}] = 0
\end{equation}
in simply connected domains $\ \Lambda\subset{\mathbf R^8}$
$(\Lambda\cap\mathbf R=\emptyset),\ $ 
with the  $C^2$-boundary $\partial{\Lambda}$, 
to within arbitrary constant. 

 \section{Conclusions}
  
New axially symmetric system $(A_n)$ takes up an intermediate place 
between spherically symmetric system investigated 
by Brackx, Delange and Sommen \cite{BrDeSo:1982} and 
an asymmetric system $(H_n)$ investigated by Leutwiler \cite{Leut:CV17}.

Our approach allows to demonstrate, in particular, some transitions 
between lower and higher dimensions for quaternionic and 
octonionic generalizations of classical holomorphic functions. 

How many various generalizations of the Cauchy-Riemann system 
in $\mathbf R^8$ having solutions in the form of functions of 
the octonionic variable exist ? This is an open question. 


{\bf Acknowledgement.\ } The author would like to thank Dr.Kisil
for technical support.

 {\bf Address}: 19 Mira av., apt. 225 Fryazino 
 
                Moscow region 141190 Russia 

 {\bf e-mail}:  bryukhov@mail.ru, bryukhov@pochtamt.ru 


\small
\begin{thebibliography}{}

 \bibitem {Bit:BoEl} Bitsadze A.~V., 
\emph{"Boundary Value Problems for Second Order Elliptic Equations"},
 Amsterdam, North-Holland, 1968 
(Russian Original: Moscow, Nauka, 1966).

 \bibitem {BrDeSo:1982} Brackx F., Delange R. and F.~Sommen,
\emph{"Clifford Analysis"}, Research Notes in Math. \textbf{76},
Boston-London-Melbourne, Pitman, 1982.



 \bibitem {GilTru} Gilbarg D. and N.~S. Trudinger,
\emph{"Elliptic Partial Differential Equations of Second Order"},
 2nd Edition, Berlin, Springer, 1983.


 \bibitem {HeLe:1996}
 Hempfling~Th. and H.~Leutwiler, 
\emph{Modified quaternionic analysis in $\mathbf R^4$},
"Proc. of Clifford Algebras and Their Appl. in Math. Physics (Aachen, 1996)".
Fundamental Theories of Physics \textbf{94}, 227--237,
Dordrecht, Kluwer Academic Publishers, 1998.





 \bibitem {Landis} Landis E.~M., 
 \emph{"Second order equations of elliptic and parabolic types"},
 Providence, RI, American Mathematical Society, 1998
 (Russian Original: Moscow, Nauka, 1971).

 \bibitem {LavSh} Lavrentiev M.~A. and B.~V. Shabat,
\emph{"Methods of the theory of functions in a complex variable"},
5th Edition, Moscow, Nauka, 1987 (in Russian).

 \bibitem {Leut:CV17} Leutwiler~H., 
\emph{Modified Clifford analysis},
Complex Variables Theory Appl.
\textbf{17}, (1992) 153--171.

 \bibitem {Leut:CV20} Leutwiler~H.,
\emph{Modified quaternionic analysis in $\mathbf R^3$},
Complex Variables Theory Appl.
\textbf{20}, (1992) 19--51.

 \bibitem {Leut:FM7} Leutwiler~H., 
\emph{More on modified quaternionic analysis in $\mathbf R^3$},
Forum Math. \textbf{7}, (1995) 279--305.

 \bibitem {Leut:EM14} Leutwiler~H., 
\emph{Rudiments of function theory in $\mathbf R^3$},
Expositiones Math. \textbf{14}, (1996) 97--123.

 \bibitem {LiKaPe:2001} Xingmin Li, Zhao Kai and Lizhohg Peng, 
\emph{The Laurent series on the octonions}, 
"Int. Conf. on Clifford Analysis, Its Appl. and Related Topics. Beijing,
August 1-6, 2000".
Adv. Appl. Clifford Algebras \textbf{11(S2)}, (2001) 205--211.

 \bibitem {LiPe:2001} Xingmin Li and Lizhohg Peng, 
\emph{Taylor series and orthogonality of the octonion analytic functions},
Acta of Math. Scientia \textbf{3}, (2001).



 \bibitem {PolBr} van der Pol B. and H. Bremmer, 
\emph{"Operational Calculus based on the Two-Sided Laplace Integra"}, 
2nd Edition, Cambridge, Cambridge Univ. Press, 1955.


 \bibitem {ScTiTo:1997} Scheiher~K., Tichy~R.F. and K.W. Tomantschger, 
\emph{Elementary inequalities in hypercomplex numbers},
Anz. Osterreich. Akad. Wiss. Math.-Natur. Kl. \textbf{134}, (1997) 3--10.



 \bibitem {Titch:R} Titchmarsh E.~C., 
\emph{"The Theory of the Riemann Zeta-Function"}, 
2nd Edition, Oxford, Oxford Univ. Press, 1986.


 \bibitem {WW} Whittaker E.~T. and G.~N. Watson, 
\emph{"A Course of Modern Analysis"}, 
4th Edition, Cambridge, Cambridge Univ. Press, 1996.

 \end{thebibliography}
\end{document}